\documentclass[11pt]{amsart}
\usepackage{amssymb}
\usepackage{amsbsy}
\usepackage{epsfig}
\usepackage{graphicx}
\usepackage{tikz}
\usepackage{url}

\AtBeginDocument{%
   \def\MR#1{}
}

\DeclareMathOperator*{\diag}{diag}

\begin{document}
\pagestyle{myheadings}
\title{A conjecture on Laplacian eigenvalues of trees}
\author{David P. Jacobs}
\address{ School of Computing, Clemson University Clemson, SC 29634 USA}
\email{\tt dpj@clemson.edu}
\author{Vilmar Trevisan}
\address{Instituto de Matem\'atica, UFRGS,  91509--900 Porto Alegre, RS, Brazil}
\email{\tt trevisan@mat.ufrgs.br}
\pdfpagewidth 8.5 in
\pdfpageheight 11 in

\def\floor#1{\left\lfloor{#1}\right\rfloor}

\newcommand{\twobytwo}[4]{
\left [
    \begin{array}{lr}
       #1 & #2  \\
       #3 & #4  \\
    \end{array}
\right ]
}
\newcommand{\Real}{{\mathbb R}}
\newcommand{\theratio}[2]{ \lceil \frac{#1}{#2} \rceil }
\newcommand{\displayratio}[2]{ \left \lceil \frac{#1}{#2} \right \rceil }
\newcommand{\avgdeg}[1]{2 - \frac{2}{#1} }
\newcommand{\floorratio}[2]{ \left\lfloor \frac{#1}{#2} \right\rfloor }
\newcommand{\Prf}{{\bf Proof: }}
\newcommand{\PrfSketch}{{\bf Proof (Sketch): }}
\newcommand{\boldQ}{\mbox{\bf Q}}
\newcommand{\boldR}{\mbox{\bf R}}
\newcommand{\boldZ}{\mbox{\bf Z}}
\newcommand{\boldc}{\mbox{\bf c}}
\newcommand{\sign}{\mbox{sign}}
\newtheorem{Thr}{Theorem}
\newtheorem{Pro}{Proposition}
\newtheorem{Con}{Conjecture}
\newtheorem{Cor}{Corollary}
\newtheorem{Lem}{Lemma}
\newtheorem{Fac}{Fact}
\newtheorem{Ex}{Example}
\newtheorem{Def}{Definition}
\newtheorem{Prob}{Problem}
\def\floor#1{\left\lfloor{#1}\right\rfloor}

\newenvironment{my_enumerate}{
\begin{enumerate}
  \setlength{\baselineskip}{14pt}
  \setlength{\parskip}{0pt}
  \setlength{\parsep}{0pt}}{\end{enumerate}
}
\newenvironment{my_description}{
\begin{description}
  \setlength{\baselineskip}{14pt}
  \setlength{\parskip}{0pt}
  \setlength{\parsep}{0pt}}{\end{description}
}

\begin{abstract}
Motivated by classic tree algorithms, in 1995 we designed a bottom-up
$O(n)$ algorithm to compute the determinant of a tree's adjacency matrix $A$.
In 2010 an $O(n)$ algorithm was found for constructing
a diagonal matrix congruent to $A + xI_n$, $x \in \Real$,
enabling one to easily count the number of eigenvalues in any interval.
A variation of the algorithm allows Laplacian eigenvalues in trees to be counted.
We conjecture that for any tree $T$ of order $n \geq 2$,
at least half of its
Laplacian eigenvalues are less than $\bar{d} = \avgdeg{n}$,
its average vertex degree.
\end{abstract}

\thanks{Work supported by CNPq - Grant 400122/2014-6, Brazil}

\maketitle

\section{Background}
\label{sec-background}
The goal of this paper is to describe a fascinating conjecture in spectral graph theory
involving Laplacian eigenvalues in trees,  as well as some eigenvalue location algorithms that
we designed.  But before describing these things, we want to give some background,
recalling some early influences that led us on this path.

It is worth noting that algorithms for trees have been around a long time.
Recalling that a set of vertices in a graph is {\em independent} if no two members are adjacent,
in 1966 Daykin and Ng \cite{DaykinNg66}
gave the first algorithm for computing
$\beta_0$,
the size of a
largest independent set in a tree $T$.
A vertex set $S$ is {\em dominating} if every vertex $v \not\in S$
is adjacent to some member of $S$.
In 1975, Cockayne, Goodman and Hedetniemi \cite{Cockayne75}
published the first algorithm for computing
$\gamma$, the size of a smallest dominating set in a tree.
See \cite{Hedet2006} for more background on tree algorithms.
Many tree algorithms are bottom-up, carrying information to the root.

Beginning in the 1980's and for about thirty years,
Steve Hedetniemi
led a weekly algorithms
seminar at Clemson University.  In the seminar's early days, a major theme was finding
linear-time algorithms for NP-complete problems, when restricted to a graph class such as trees.
One of the frequent presenters at the seminar was Tom Wimer,
one of Hedetniemi's Ph.D. students,
who had developed a methodology for linear-time algorithms \cite{Wimer}.

In 1995 Jacobs took a sabbatical to work with Trevisan's computer algebra group
in the Mathematics Institute at
Universidade Federal do Rio Grande do Sul
(UFRGS) in Brazil.
UFRGS is the top research university in the region,
located in Porto Alegre, a city of about 1.4 million and the capital
of Rio Grande do Sul, the southernmost state.
It was during that year that our collaboration began.
While we spent many long hours in front of the blackboard,
we also attended Gr\^emio games,
ate churrasco, and saw some beautiful beaches.

A problem we worked on was finding a fast algorithm to compute the determinant
of a tree's adjacency matrix.
Recall that the {\em adjacency matrix} $A$ of a graph $G =(V,E)$
is the matrix indexed by $V \times V$, whose $(u,v)$
entry is $1$ if $\{u,v\} \in E$, and $0$ otherwise.
While the {\em determinant} of a matrix is often defined
inductively with expansion by co-factors \cite{MarcusMinc},
this is a bad way to compute it.
An $O(n^3)$ method is to apply certain operations,
reducing it to a triangular matrix, whose determinant is the product of the diagonal.
Adding multiples of a row (column) to another row (column)
does not alter the determinant.
Interchanging two rows or columns simply changes the determinant's sign.
For arbitrary $n \times n$ matrices the complexity of computing determinants
is known to be that of matrix multiplication \cite{Burgisser}, which
can be done presently in $O(n^{2.376})$.

One of our nicest achievements that year was an $O(n)$ algorithm \cite{FHJT96}
for computing the determinant of
the adjacency matrix of a tree $T$.
The algorithm works bottom-up like many graph algorithms,
and is similar to the  tree algorithms that will be shown later in this paper.
To achieve linear time (and linear space) the algorithm could not store the entire matrix.
Rather, the algorithm operated on the tree itself, beginning at the leaves, and performing
calculations as it worked its way to the root.
The numbers stored in the nodes of the tree were values
that would have been diagonal elements, had we used ordinary
row and column operations.
In essence, the adjacency matrix was being transformed
to a diagonal matrix using elementary row and column operations.
At the end of this diagonalization process, the determinant is
to within a sign, the product of the values in the tree nodes.

After we saw how to do this, we sought solutions to other algebraic problems involving the trees.
These included finding the LUP decomposition \cite{JT1999} of the adjacency matrix,
computing the characteristic polynomial \cite{JT1998,JT2005} of the adjacency matrix,
and computing the inverse of a tree's incidence matrix \cite{Jacobs2008} with Machado and Pereira.
While most of our work has been algebraic, we also wrote a beautiful number theoretic paper
characterizing so-called Chebyshev numbers \cite{JTR2008}, and discovered a primality testing
algorithm based on Chebyshev polynomials \cite{JTR2008a}, both co-written with Mohammad Rayes.

In the next section we will give some fundamental ideas in spectral graph theory,
and highlight some classic results.
In Section~\ref{sec-adjalgorithm} we will give our $O(n)$ tree algorithm
for obtaining a diagonal matrix congruent to $A + xI_n$, $x \in \Real$,
and explain its use in finding eigenvalues of trees.
The Laplacian matrix and the algorithm's Laplacian analog are
given in Section~\ref{sec-laplacealgorithm},
along with some classic theorems involving Laplacian eigenvalues.
Finally, in Section~\ref{sec-conjecture}
we will discuss our conjecture on the distribution
of Laplacian eigenvalues for trees.

\section{Spectral Graph Theory - SGT}
\label{sec-specgraphthry}
The main goal of {\em spectral graph theory} - SGT - sometimes called
{\em algebraic graph theory} is to study structural properties of a graph based on its spectrum,
that is, the eigenvalues of the matrices associated with it.
The first papers appeared in connection with quantum chemistry
(the 1931 paper  by H\"uckel \cite{Huc31} is considered the first published work on the subject).
Later, SGT was popularized by D. Cvetkovi\'{c}'s 1971 PhD thesis  \cite{Cve71}.
Nowadays, SGT is a fully developed field that has ingredients from both graph theory and
linear algebra, but has its own characteristic features.
The web page \cite{SGT}
provides useful and up to date information on
many aspects of the theory of graph spectra for researchers interested in
the field and its applications.

In this section, we will primarily deal with the adjacency matrix, but there are many other
matrices associated with a graph $G$. Our goal is to introduce a few well-known results that
are representative of SGT and that relate spectral parameters
to classical graph theory.

Recall that the \emph{eigenvalues} of a square matrix $M$ are the
roots of the \emph{characteristic polynomial} $\det(\lambda I - M)$.
The \emph{spectrum} of $M$ is the (multi)set determined by the roots of the
characteristic polynomial of $M$, together with their multiplicities.
Since the $n \times n$ adjacency matrix $A$ of $G$ is real and symmetric, its
eigenvalues are real and we order them as
$$
\lambda_1 \ge \lambda_2 \ge \cdots \ge \lambda_n .
$$
We observe that the characteristic polynomial of a graph $G$ is invariant under the vertex ordering,
hence the spectrum of $G$ is independent of the ordering used for its vertices.

\begin{figure}[t]
\begin{minipage}[t]{0.8\textwidth}
  \begin{tikzpicture}
    [scale=1,auto=left,every node/.style={circle,minimum width={.25in}, scale=0.7}]
   \node[draw,circle,fill=white] (a) at (0,0) {};
   \node[draw,circle,fill=white] (b) at (-1.5,0) {};
   \node[draw,circle,fill=white] (c) at (1.5,0) {};
   \node[draw,circle,fill=white] (d) at (0,1.5) {};
   \node[draw,circle,fill=white] (e) at (0,-1.5) {};
   \path
       (a) edge [thick] node[left]{} (b)
       (a) edge [thick] node[left]{} (c)
       (a) edge [thick] node[right]{} (d)
       (a) edge [thick] node[right]{} (e);
 \end{tikzpicture}
  \begin{tikzpicture}
    [scale=1,auto=left,every node/.style={circle,minimum width={.35in}, scale=0.7}]
    \node[draw,circle,fill=none,opacity=0.0] (f) at (7.5,0) {};
  \end{tikzpicture}
  \begin{tikzpicture}
    [scale=1,auto=left,every node/.style={circle,minimum width={.35in}, scale=0.7}]
    \node[draw,circle,fill=none,opacity=0.0] (f) at (7.5,0) {};
  \end{tikzpicture}
  \begin{tikzpicture}
    [scale=1,auto=left,every node/.style={circle,minimum width={.25in}, scale=0.7}]
    \node[draw,circle,fill=white] (f) at (6.5,0) { };
    \node[draw,circle,fill=white] (g) at (5.0,1.5) { };
    \node[draw,circle,fill=white] (h) at (8.0,1.5) { };
    \node[draw,circle,fill=white] (i) at (5.0,-1.5) { };
    \node[draw,circle,fill=white] (j) at (8.0,-1.5) { };
    \path
        (g) edge [thick] node[right]{} (h)
        (i) edge [thick] node[right]{} (j)
        (g) edge [thick] node[right]{} (i)
        (j) edge [thick] node[right]{} (h);
 \end{tikzpicture}
\end{minipage}
\caption{\label{coespec} Nonisomorphic graphs with the same spectrum.}
\end{figure}
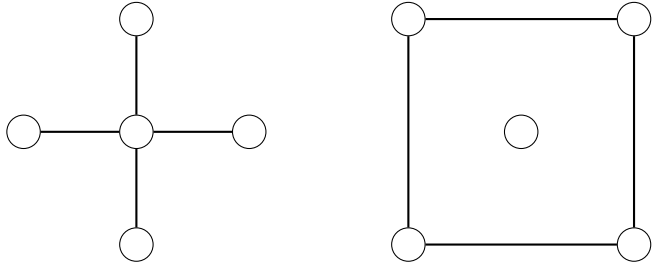

Given a graph $G$, its \emph{spectrum} is the spectrum of its adjacency matrix.
A graph $G$ is said to be \emph{determined by its spectrum} if any
other graph, nonisomorphic  to $G$, has a different spectrum.
It is well-known that graphs in general are not determined by their spectra.
The example in Figure \ref{coespec} shows two non-isomorphic graphs of order $n = 5$,
both having the characteristic polynomial  $p(\lambda)=\lambda^5 - 4\lambda^3$ and
spectrum $\{2,0,0,0,-2\}$. A typical  problem in SGT is to determine classes of
graphs that are determined by their spectra (see, for example, \cite{van2003}).
On the other hand, there exist parameters that are determined by the spectrum of the graph $G$.
Besides the (obvious) fact that the number of vertices equals the degree
of the characteristic polynomial,
the number of triangles of a
graph is also determined by the spectrum.
If
$$
p_G(\lambda) = \lambda^n + a_1 \lambda^{n - 1} + a_2 \lambda^{n - 2} + \cdots + a_{n - 1} \lambda + a_n
$$
is the characteristic polynomial of an order $n$ graph $G$ with $m$ edges, then
 \begin{itemize}
\samepage
\item[(i)] $a_1 = 0$;
\item[(ii)] $a_2 = - m$;
\item[(iii)] $a_3 = - 2 t$, where  $t$ is the number of triangles in $G$.
\end{itemize}
This can be proven by analyzing the relationship between the coefficients of the characteristic
polynomial and the principal minors of the adjacency matrix.
Moreover, using the fact that the number of walks of length $k$ from
vertex $v_i$ to vertex $v_j$ is the $i,j$ entry of $A^k$, and that the trace is
the sum of the eigenvalues of a matrix, we see that
\begin{itemize}
\samepage
 \item[(i)] $\sum_{i=1}^n {\lambda{_i}}^2  = 2m$;
 \item[(ii)] $\sum_{i=1}^n {\lambda{_i}}^3  = 6t$.
\end{itemize}
However, the number of cycles of length 4 or more,
is not determined the by a graph's spectrum as
the graphs in Figure \ref{coespec} show.

The {\em distance} between two vertices $u$ and $v$ is
the number of edges in a shortest path between them, and the graph's {\em diameter}
is the greatest distance between any two vertices.
\begin{Thr}
Let $G$ be a connected graph with diameter $d$.
The number of distinct eigenvalues of $G$ is at least $d+1$.
\end{Thr}

\noindent
\textbf{Proof:}
Let $\lambda_1, \lambda_2, \ldots, \lambda_t$ be the distinct eigenvalues of $G$.
Since $A$ is real and symmetric, the minimal polynomial of $A$ has degree $t$ and
$$
(A - \lambda_1 I)(A -\lambda_2 I)\dots (A - \lambda_t I) = 0.
$$
Hence, $A^t$ is a linear combination of $I, A, \ldots, A^{t-1}$.
Let  $u,v$ be two vertices whose distance is $d$.
Then the entry $[A^i]_{u,v} = 0$
for all $i=0,\ldots, d-1$ and the entry of $[A^d]_{u,v} > 0$.
If $t \leq d$ this would imply a contradiction $\Box$.

By the Perron-Frobenius theory, it is known that
the largest eigenvalue $\lambda_1$
is positive,
simple (has multiplicity $1$),
is the spectral radius $\rho(A) = \max_i |\lambda_i|$,
and has an eigenvector with all positive components.
The largest eigenvalue
is called the graph's \emph{index}, and is an
important spectral parameter as the next results illustrate.

An early result
is the following remarkable
observation due to Wilf \cite{Wil67}.
Recall that the {\em chromatic number} of a graph is the smallest number of colors
needed to properly color its vertices.
\begin{Thr} Let $G$ be an $n$-vertex connected graph with
chromatic number $\chi$ and index $\lambda_1$. Then $$\chi \leq 1 + \lambda_1,$$
with equality iff $G$ is the complete graph or an odd cycle.
\end{Thr}

Using the fact that $Tr(A) = \sum_i^n \lambda_i = 0$ and that
$Tr(A^2) = \sum_i^n \lambda_i^2 = 2m$, where $m$ is the number of edges
of $G$, Wilf proves the following upper bound.
\begin{Cor}
If $G$ has $m$ edges and $n$ vertices, then
$$
\chi \leq \left(2(1-\frac{1}{n})m \right)^{\frac{1}{2}} +1.
$$
\end{Cor}

A lower bound for the chromatic number of a graph, in terms of spectral
parameters, given by Hoffman \cite{Hoff70} is
$$
\chi \geq 1+ \frac{\lambda_1}{-\lambda_n}.
$$
Nikiforov \cite{Nik2007} improved this bound
using $\mu_1$, the
largest Laplacian eigenvalue, (defined in Section \ref{sec-laplacealgorithm}) showing
$$
\chi \geq 1+ \frac{\lambda_1}{\mu_1 - \lambda_n}.
$$

Important relations between the index and the hamiltonicity
of a graph have been studied by several authors.
To state the most current results,
we write $K_{n-1} + v$ for the complete graph on $n - 1$ vertices
together with an isolated vertex $v$, and $K_{n-1} + e$ for the complete graph on
$n - 1$ vertices with a pendent edge $e$.
The following results are due to Fiedler and Nikiforov \cite{Fie2010} and
represent improvements on previous bounds, whose references may be found there.

\begin{Thr}
Let $G$ be a graph of order $n$ and index $\lambda_1$. If
$\lambda_1 \geq n-2$,
then $G$ contains a hamiltonian path unless $G = K_{n-1} + v$.
If the inequality is strict, then $G$ contains a hamiltonian cycle unless $G = K_{n-1} + e$.
\end{Thr}

\begin{Thr}
Let $G$ be a graph of order $n$ and let
$\bar{\lambda_1}$
denote
the index of the complement of $G$.
If
$$
\bar{\lambda_1} \leq \sqrt{n - 1},
$$
then $G$ contains a hamiltonian path unless $G = K_{n-1} + v$.
If
$$
\bar{\lambda_1} \leq \sqrt{n - 2},
$$
then $G$ contains a hamiltonian cycle unless $G = K_{n-1} + e$.
\end{Thr}

We finish this section by mentioning an application of SGT to network problems.
Locating influential vertices in a network has many applications such as vertex partitions,
social networks and disease control, and the literature is really abundant on the subject.
We want to mention that the index of the network may be used to measure the importance of its vertices.
Bonacich \cite{Bon87} defined the \emph{eigenvector centrality measure} of
a vertex $v_i$ as the $i$-th component of the eigenvector
$\textbf{v}$
associated to the index $\lambda_1$. The intuitive idea is that larger values
in $\textbf{v}$ indicate more important vertices.

\section{New kind of algorithm}
\label{sec-adjalgorithm}
In 2010 we met for several weeks in Clemson to develop a new kind of algorithm.
The algorithm that we designed takes, as input,
a tree $T$ and $x \in \Real$, and
constructs a diagonal matrix $D$ congruent to $A + x I_n$,
where $A$ is the adjacency matrix of $T$.
Two matrices $R$ and $S$ are {\em congruent}
if there exists a nonsingular matrix $P$ such that $S = P^T R P$.
Or equivalently, $R$ can be transformed to $S$
by applying the {\em same} elementary row and column operations.
For example,
$$
R = \twobytwo{2}{1}{1}{2},  S = \twobytwo{2}{0}{0}{\frac{3}{2}}
$$
are congruent because $S = P^T R P$
where
$$
P = \twobytwo{1}{-\frac{1}{2}}{0}{1}
$$
or equivalently, one may obtain
$S$ from $R$ by subtracting half the first row from the second,
and then applying the same column operation.

Like the determinant algorithm that inspired it, the algorithm
operates directly on the tree, hiding row and column operations.
It begins by initializing all vertices with $x$.
Then working bottom-up toward the root, each non-leaf $v$ is processed.
Whenever possible, the algorithm simply subtracts from $v$'s value,
the {\em reciprocal} values of $v$'s children.
\begin{centering}
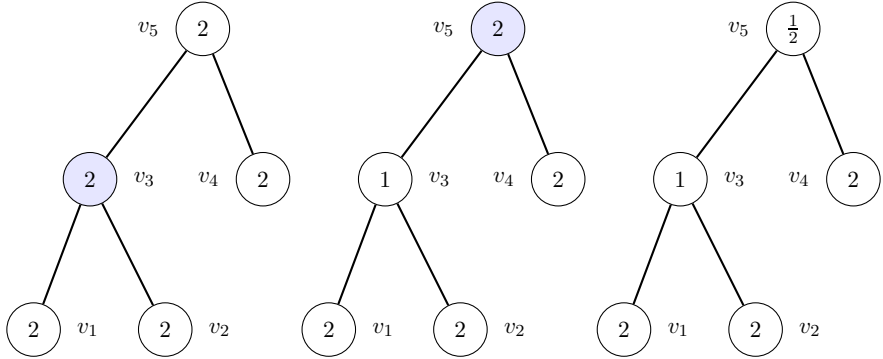
\begin{figure}[t]
\begin{minipage}[t]{0.3\textwidth}
\begin{tikzpicture}
  [scale=1,auto=left,every node/.style={circle,minimum width={.35in}, scale=0.8}]
  \node[draw,circle,fill=white, label=left:$v_5$] (a) at (0,2) { 2 };
  \node[draw,circle,fill=blue!10, label=right:$v_3$] (b) at (-1.5,0) { 2 };
  \node[draw,circle,fill=white, label=left:$v_4$] (c) at (.8,0) { 2 };
  \node[draw,circle,fill=white, label=right:$v_1$] (d) at (-2.25,-2) { 2 };
  \node[draw,circle,fill=white, label=right:$v_2$] (e) at (-.5, -2) { 2 };
  \path 
      (a) edge [thick] node[left]{} (b)
      (a) edge [thick] node[left]{} (c)
      (b) edge [thick] node[right]{} (d)
      (b) edge [thick] node[right]{} (e);
\end{tikzpicture}
\end{minipage}
\begin{minipage}[t]{0.3\textwidth}
\begin{tikzpicture}
  [scale=1,auto=left,every node/.style={circle,minimum width={.35in}, scale=0.8}]
  \node[draw,circle,fill=blue!10, label=left:$v_5$] (a) at (0,2) { 2 };
  \node[draw,circle,fill=white, label=right:$v_3$] (b) at (-1.5,0) { 1 };
  \node[draw,circle,fill=white, label=left:$v_4$] (c) at (.8,0) { 2 };
  \node[draw,circle,fill=white, label=right:$v_1$] (d) at (-2.25,-2) { 2 };
  \node[draw,circle,fill=white, label=right:$v_2$] (e) at (-.5, -2) { 2 };
  \path 
      (a) edge [thick] node[left]{} (b)
      (a) edge [thick] node[left]{} (c)
      (b) edge [thick] node[right]{} (d)
      (b) edge [thick] node[right]{} (e);
\end{tikzpicture}
\end{minipage}
\begin{minipage}[t]{0.3\textwidth}
\begin{tikzpicture}
  [scale=1,auto=left,every node/.style={circle,minimum width={.35in}, scale=0.8}]
  \node[draw,circle,fill=white, label=left:$v_5$] (a) at (0,2) { $\frac{1}{2}$ };
  \node[draw,circle,fill=white, label=right:$v_3$] (b) at (-1.5,0) { 1 };
  \node[draw,circle,fill=white, label=left:$v_4$] (c) at (.8,0) { 2 };
  \node[draw,circle,fill=white, label=right:$v_1$] (d) at (-2.25,-2) { 2 };
  \node[draw,circle,fill=white, label=right:$v_2$] (e) at (-.5, -2) { 2 };
  \path 
      (a) edge [thick] node[left]{} (b)
      (a) edge [thick] node[left]{} (c)
      (b) edge [thick] node[right]{} (d)
      (b) edge [thick] node[right]{} (e);
\end{tikzpicture}
\end{minipage}
\caption{\label{examplefig1} Diagonalization with $x = 2$ }
\end{figure}
\end{centering}

Figure~\ref{examplefig1} shows a tree rooted at $v_5$, and depicts the execution of
the algorithm when $x = 2$.
The leftmost diagram shows the tree after initialization,
and $v_3$ is shaded to indicate it will be the first nonleaf to be processed.
The center diagram shows the new value $2 - \frac{1}{2} - \frac{1}{2}$ in $v_3$,
and indicates that $v_5$ will be the next vertex to be processed.
The rightmost diagram gives the updated value $2 - 1 - \frac{1}{2}$ in $v_5$,
and shows $\diag(D) = (2,2,1,2,\frac{1}{2})$.

\begin{centering}
\begin{figure}[h]
\begin{minipage}[t]{0.3\textwidth}
\begin{tikzpicture}
  [scale=1,auto=left,every node/.style={circle,minimum width={.35in}, scale=0.8}]
  \node[draw,circle,fill=white, label=left:$v_5$] (a) at (0,2) { 0 };
  \node[draw,circle,fill=blue!10, label=right:$v_3$] (b) at (-1.5,0) { 0 };
  \node[draw,circle,fill=white, label=left:$v_4$] (c) at (.8,0) { 0 };
  \node[draw,circle,fill=white, label=right:$v_1$] (d) at (-2.25,-2) { 0 };
  \node[draw,circle,fill=white, label=right:$v_2$] (e) at (-.5, -2) { 0 };
  \path 
      (a) edge [thick] node[left]{} (b)
      (a) edge [thick] node[left]{} (c)
      (b) edge [thick] node[right]{} (d)
      (b) edge [thick] node[right]{} (e);
\end{tikzpicture}
\end{minipage}
\begin{minipage}[t]{0.3\textwidth}
\begin{tikzpicture}
  [scale=1,auto=left,every node/.style={circle,minimum width={.35in}, scale=0.8}]
  \node[draw,circle,fill=blue!10, label=left:$v_5$] (a) at (0,2) { 0 };
  \node[draw,circle,fill=white, label=right:$v_3$] (b) at (-1.5,0) {$-\frac{1}{2}$ };
  \node[draw,circle,fill=white, label=left:$v_4$] (c) at (.8,0) { 0 };
  \node[draw,circle,fill=white, label=right:$v_1$] (d) at (-2.25,-2) { 2 };
  \node[draw,circle,fill=white, label=right:$v_2$] (e) at (-.5, -2) { 0 };
  \path 
      (a) edge [dashed] node[left]{} (b)
      (a) edge [thick] node[left]{} (c)
      (b) edge [thick] node[right]{} (d)
      (b) edge [thick] node[right]{} (e);
\end{tikzpicture}
\end{minipage}
\begin{minipage}[t]{0.3\textwidth}
\begin{tikzpicture}
  [scale=1,auto=left,every node/.style={circle,minimum width={.35in}, scale=0.8}]
  \node[draw,circle,fill=white, label=left:$v_5$] (a) at (0,2) { $-\frac{1}{2}$ };
  \node[draw,circle,fill=white, label=right:$v_3$] (b) at (-1.5,0) { $-\frac{1}{2}$ };
  \node[draw,circle,fill=white, label=left:$v_4$] (c) at (.8,0) { 2 };
  \node[draw,circle,fill=white, label=right:$v_1$] (d) at (-2.25,-2) { 2 };
  \node[draw,circle,fill=white, label=right:$v_2$] (e) at (-.5, -2) { 0 };
  \path 
      (a) edge [dashed] node[left]{} (b)
      (a) edge [thick] node[left]{} (c)
      (b) edge [thick] node[right]{} (d)
      (b) edge [thick] node[right]{} (e);
\end{tikzpicture}
\end{minipage}
\caption{\label{examplefig2} Diagonalization with $x = 0$ }
\end{figure}
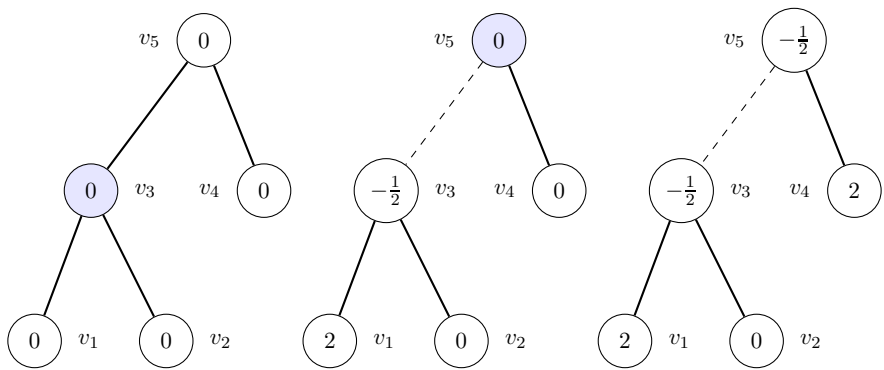
\end{centering}

Reciprocals require numbers be nonzero.
Should a vertex have a child whose value is zero,
it chooses one such child,
assigns itself $-\frac{1}{2}$ and the child $2$,
and breaks the edge with its parent (unless it is the root).
All other children are unaffected.
Figure ~\ref{examplefig2} depicts the execution of the algorithm when $x = 0$.
In the rightmost diagram, we see $\diag(D) = (2,0,-\frac{1}{2},2, -\frac{1}{2})$.

The algorithm is shown in Figure~\ref{alg-adj}.
It assumes vertices are numbered so that if $v_i$ is a child of $v_j$ then $i < j$.
A correctness proof and some applications can be found it \cite{JT2011}.
Clearly, it makes $O(n)$ arithmetic operations.
\begin{figure}[h]
\input{fig-alg-adj.tex}
\caption{\label{alg-adj} $A + x I_n$ diagonalization}
\end{figure}

There is more to this story.
We can use the algorithm Diagonalize
to study and locate eigenvalues of trees.
Using a classic result
known as Sylvester's Law of Inertia,
one can show (see \cite{JT2011}):
\begin{Thr}
\label{thr-inertia}
Let $D$ be the diagonal matrix returned by $\text{Diagonalize}(T,-x)$.
\begin{enumerate}
\samepage
\item The number of positive entries of $D$
is the number of eigenvalues of $T$ greater than $x$.
\item The number of negative entries of $D$
is the number of eigenvalues of $T$ less than $x$.
\item The number of zero entries in $\diag(D)$ is the multiplicity of $x$.
\end{enumerate}
\end{Thr}

Applying Theorem~\ref{thr-inertia} to the diagonalization
in Figure~\ref{examplefig2}, we conclude that the tree has two
positive eigenvalues, two negative eigenvalues and that zero has multiplicity one.
From Figure~\ref{examplefig1}, we conclude that all its eigenvalues
are greater than $-2$.
The tree's spectrum is actually $0, \pm \sqrt{ 2 \pm \sqrt{2} }$.

\section{Laplacian eigenvalues}
\label{sec-laplacealgorithm}
Let $G = (V,E)$, be an undirected graph with vertex set
$V = \{v_1, \dots, v_n\}$, and adjacency matrix $A$.
The {\em Laplacian matrix} of $G$ is given by $L = D - A$, where $D = [d_{ij}]$
is the diagonal matrix in which $d_{ii} = \deg(v_i)$, the degree of $v_i$.
The {\em Laplacian spectrum} of $G$ is the set of eigenvalues
of $L$, which we will denote as:
$$
\mu_1 \geq \mu_2 \geq \ldots \geq \mu_n .
$$
It can be shown that for any graph,
$\mu_n = 0$ and $\mu_1 \leq n$,
and so Laplacian eigenvalues are nonnegative.
Moreover, $\sum_{i=1}^n \mu_i = 2m$.
See \cite{Merr94,mohar,mohar92} for surveys on many interesting properties
of the Laplacian spectrum.

As we saw in Section~\ref{sec-specgraphthry},
the eigenvalues of a graph's adjacency matrix
can be used to infer properties about the graph.
Similarly, a graph's Laplacian eigenvalues also give information
about the graph.

The following beautiful result,
sometimes called the Matrix Tree Theorem,
is attributed to the nineteenth century mathematician G. Kirchhoff
\cite{Kirchhoff}.
More information can be found in \cite{ChaikenKleitman},
and a proof can be found in \cite{BrouwerHaemers}.

\begin{Thr}
\label{MatrixTreeThr}
The number of spanning trees in a graph $G$ is given by
$$
\frac{ \mu_1 \mu_2 \ldots \mu_{n-1} }{n} .
$$
\end{Thr}

The second smallest Laplacian eigenvalue $\mu_{n-1}$ plays a fundamental role
in spectral graph theory, and is called the graph's \emph{algebraic connectivity}.
We can see that if $\mu_{n-1}=0$, then the product in Theorem~\ref{MatrixTreeThr} is $0$,
and so the graph is disconnected.
More generally, it can be shown that the multiplicity of $0$ as a Laplacian eigenvalue
is the number of connected components of a graph $G$ and, in particular, $G$ is connected iff  $\mu_{n-1}>0$.
Much of our understanding of this spectral parameter is due Miroslav Fiedler,
whose pioneering work \cite{Fie73,Fie75} is described by Nikiforov \cite{Nik13} as a mathematical gold strike.
Fiedler studied the eigenvector associated with $\mu_{n-1}$, now called the {\em Fiedler vector},
and his findings were fundamental in understanding properties of the graph.
The papers \cite{Abr07,Abr14} contain more applications of algebraic connectivity.

There are several results, similar in spirit to Wilf's theorem,
that relate the distribution of the Laplacian spectrum
to $\gamma$, the domination number of a graph.
Lu, Liu and Tian \cite{LLT2005}
obtained a bound for $\mu_1$,
which was improved by Nikiforov \cite{Nik2007a} as follows:
\begin{Thr}
If a graph has order $n \geq 2$ and domination number $\gamma$, then
$$
\mu_1 \geq \displayratio{n}{\gamma} .
$$
\end{Thr}

If $I$ is a real interval, $m_G(I)$ denotes the number of Laplacian
eigenvalues of $G$ in $I$.
Quite recently Hedetniemi, Jacobs and Trevisan obtained in \cite{HJT2015}
\begin{Thr}
\label{thr-gamma}
If $G$ has domination number $\gamma$ then $m_G [0,1) \leq \gamma$.
\end{Thr}

One might ask if
the algorithm of Figure~\ref{alg-adj}
has a Laplacian analog.
The answer is yes.
Not long after it appeared,
it was soon discovered \cite{Braga2013} that
it can be modified to perform diagonalization on $L + xI_n$.
The only change that is required is the initialization.
Rather than initializing each vertex $v$ with $x$, we
initialize each vertex with $x + \deg(v)$.
\begin{figure}
\input{fig-alg-lap.tex}
\caption{\label{alg-lap} $L + x I_n$ diagonalization}
\end{figure}
The algorithm is shown in Figure~\ref{alg-lap}.
What is particularly nice is that the analog of
Theorem~\ref{thr-inertia} holds.
\begin{Thr}
\label{thr-lap-inertia}
Let $D$ be the diagonal matrix returned by $\text{LDiagonalize}(T,-x)$.
\begin{enumerate}
\samepage
\item The number of positive entries of $D$
is the number of Laplacian eigenvalues of $T$ greater than $x$.
\item The number of negative entries of $D$
is the number of Laplacian eigenvalues of $T$ less than $x$.
\item The number of zero entries in $\diag(D)$ is the multiplicity of $x$.
\end{enumerate}
\end{Thr}
We illustrate the algorithm in Figure~\ref{examplefig3}
and perform the diagonalization with $x = -1$.
Note how initial values depend on the degree of the vertices.
After initialization,
execution proceeds exactly as it does in
Algorithm {\tt Diagonalize} in Figure~\ref{alg-adj}.
The final diagonal is the same as in Figure~\ref{examplefig2},
but our conclusion is different.
In light of Theorem~\ref{thr-lap-inertia},
we know that $1$ is a Laplacian eigenvalue in $T$,
there are two Laplacian eigenvalues less than one,
and two greater than one.
More applications can be found in \cite{Braga2013}.

\section{The conjecture}
\label{sec-conjecture}
The main purpose of this paper is to support our belief that most of
the Laplacian eigenvalues of a tree are small.
In fact, we believe that at least
half of them are smaller than the average degree,
as stated precisely in
\begin{Con}
\label{TheConj}
For trees $T$ of order $n \geq 2$,
$m_T[0,\overline{d}) \geq \theratio{n}{2}$
where $\overline{d} = 2 - \frac{2}{n}$.
\end{Con}
Here
$\overline{d}$
is the average degree of $T$.
By remarks from Section~\ref{sec-laplacealgorithm}, the average degree $\overline{d}$
is also
the average Laplacian eigenvalue.
Note, however, that in general, for a graph $G$
one can have $m_G[0,\overline{d}) < \theratio{n}{2}$.
For example, the complete graph
$K_n$ has only a single Laplacian eigenvalue, namely $0$, in
$[0,\overline{d})$.
The remaining $n-1$ Laplacian eigenvalues are equal to $n$.
Another observation is that $\overline{d} = 2 - \frac{2}{n}$ can not be a Laplacian eigenvalue
of any tree with $n > 2$ vertices.
This follows because the Laplacian characteristic polynomial has integer coefficients
and, hence, any rational root must be integer.
We also note that the bound is tight, as there are trees $T$ for which
$m_T[0,\overline{d}) = \theratio{n}{2}$.
Examples of such trees are the paths $P_n$.
\begin{centering}
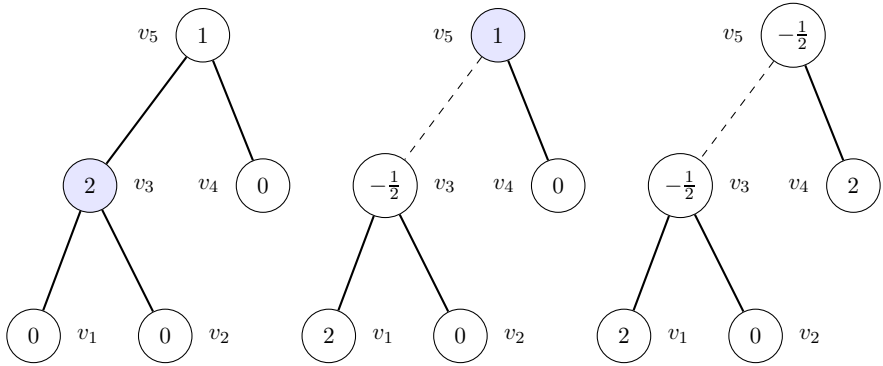
\begin{figure}[h]
\begin{minipage}[t]{0.3\textwidth}
\begin{tikzpicture}
  [scale=1,auto=left,every node/.style={circle,minimum width={.35in}, scale=0.8}]
  \node[draw,circle,fill=white, label=left:$v_5$] (a) at (0,2) { 1 };
  \node[draw,circle,fill=blue!10, label=right:$v_3$] (b) at (-1.5,0) { 2 };
  \node[draw,circle,fill=white, label=left:$v_4$] (c) at (.8,0) { 0 };
  \node[draw,circle,fill=white, label=right:$v_1$] (d) at (-2.25,-2) { 0 };
  \node[draw,circle,fill=white, label=right:$v_2$] (e) at (-.5, -2) { 0 };
  \path 
      (a) edge [thick] node[left]{} (b)
      (a) edge [thick] node[left]{} (c)
      (b) edge [thick] node[right]{} (d)
      (b) edge [thick] node[right]{} (e);
\end{tikzpicture}
\end{minipage}
\begin{minipage}[t]{0.3\textwidth}
\begin{tikzpicture}
  [scale=1,auto=left,every node/.style={circle,minimum width={.35in}, scale=0.8}]
  \node[draw,circle,fill=blue!10, label=left:$v_5$] (a) at (0,2) { 1 };
  \node[draw,circle,fill=white, label=right:$v_3$] (b) at (-1.5,0) {$-\frac{1}{2}$ };
  \node[draw,circle,fill=white, label=left:$v_4$] (c) at (.8,0) { 0 };
  \node[draw,circle,fill=white, label=right:$v_1$] (d) at (-2.25,-2) { 2 };
  \node[draw,circle,fill=white, label=right:$v_2$] (e) at (-.5, -2) { 0 };
  \path 
      (a) edge [dashed] node[left]{} (b)
      (a) edge [thick] node[left]{} (c)
      (b) edge [thick] node[right]{} (d)
      (b) edge [thick] node[right]{} (e);
\end{tikzpicture}
\end{minipage}
\begin{minipage}[t]{0.3\textwidth}
\begin{tikzpicture}
  [scale=1,auto=left,every node/.style={circle,minimum width={.35in}, scale=0.8}]
  \node[draw,circle,fill=white, label=left:$v_5$] (a) at (0,2) { $-\frac{1}{2}$ };
  \node[draw,circle,fill=white, label=right:$v_3$] (b) at (-1.5,0) { $-\frac{1}{2}$ };
  \node[draw,circle,fill=white, label=left:$v_4$] (c) at (.8,0) { 2 };
  \node[draw,circle,fill=white, label=right:$v_1$] (d) at (-2.25,-2) { 2 };
  \node[draw,circle,fill=white, label=right:$v_2$] (e) at (-.5, -2) { 0 };
  \path 
      (a) edge [dashed] node[left]{} (b)
      (a) edge [thick] node[left]{} (c)
      (b) edge [thick] node[right]{} (d)
      (b) edge [thick] node[right]{} (e);
\end{tikzpicture}
\end{minipage}
\caption{\label{examplefig3} Laplacian diagonalization with $x = -1$ }
\end{figure}
\end{centering}

At the end of this section we will state another conjecture
whose proof (or disproof) could benefit from the validity of the present conjecture.
Nevertheless, we believe
Conjecture~\ref{TheConj}
is important in itself because it helps explain
how the Laplacian eigenvalues of trees are distributed.

Conjecture~\ref{TheConj}
is known to hold for caterpillars and for paths
\cite{Braga2013}.
It is easily seen to also hold for trees whose diameter is at most 2 (stars).
As shown in \cite{Tre2011}, the conjecture holds for all trees of diameter 3 as well.

Observe that by Theorem~\ref{thr-lap-inertia},
Conjecture~\ref{TheConj}
is true if and only if for any tree $T$,
after executing
{\tt LDiagonalize}$(T,x)$
with $x=-2+\frac{2}{n}$,
at least half of the diagonal values are negative.
As an example, consider the tree $T$ of
Figure~\ref{examplefig3} whose average degree is $\frac{8}{5}$.
After executing
{\tt LDiagonalize}$(T,x)$
with $x = -\frac{8}{5}$, the final diagonal values are
$$
(-\frac{3}{5}, -\frac{3}{5}, \frac{71}{15}, -\frac{3}{5}, \frac{3041}{1065} ) .
$$
So three Laplacian eigenvalues are less than $\overline{d}$, and two are greater than $\overline{d}$.

This approach, however, may not be feasible for general trees since
the expressions resulting from
$x=-2+\frac{2}{n}$ get extremely complicated.
Nevertheless, it still may be applied to particular classes of trees,
as is illustrated by the following new result.

\begin{centering}
\begin{figure}[t]
\begin{tikzpicture}
  [scale=1,auto=left,every node/.style={circle,minimum width={.25in}, scale=0.8}]
  \node[draw,circle,fill=white, label=right:$v_0$] (a) at (0,2) {   };
  \node[draw,circle,fill=white, label=right:$v_1$] (b) at (-2.5,0) {   };
  \node[draw,circle,fill=white, label=right:$v_2 ~~~\cdots$] (c) at (-.5,0) {   };
  \node[draw,circle,fill=white, label=right:$v_r$] (g) at (2,0) {   };
  \node[draw,circle,fill=white, label=below:$s_1$] (d) at (-2.5,-2) {   };
  \node[draw,circle,fill=white, label=below:$s_2$] (e) at (-1.25, -2) {   };
  \node[draw,circle,fill=white,                  ] (f) at (-.25, -2) {   };
  \node[draw,circle,fill=white, label=below:$s_r$] (k) at (1, -2) {   };
  \node[draw,circle,fill=white,                  ] (l) at (1.75, -2) {   };
  \node[draw,circle,fill=white,                  ] (m) at (2.5, -2) {   };
  \node[draw,circle,fill=white,                  ] (n) at (3.25, -2) {   };
  \node[draw,circle,fill=white, label=above:$p$] (h) at (0,3.5) {   };
  \node[draw,circle,fill=white,                ] (i) at (-1,3.5) {   };
  \node[draw,circle,fill=white,                ] (j) at (1,3.5) {   };
  \path 
      (h) edge [thick] node[left]{} (a)
      (i) edge [thick] node[left]{} (a)
      (j) edge [thick] node[left]{} (a)
      (a) edge [thick] node[left]{} (b)
      (a) edge [thick] node[left]{} (c)
      (a) edge [thick] node[left]{} (g)
      (g) edge [thick] node[left]{} (k)
      (g) edge [thick] node[left]{} (l)
      (g) edge [thick] node[left]{} (m)
      (g) edge [thick] node[left]{} (n)
      (b) edge [thick] node[right]{} (d)
      (c) edge [thick] node[right]{} (e)
      (c) edge [thick] node[right]{} (f);
\end{tikzpicture}
\caption{\label{examplefig7} A diameter $4$ tree.}
\end{figure}
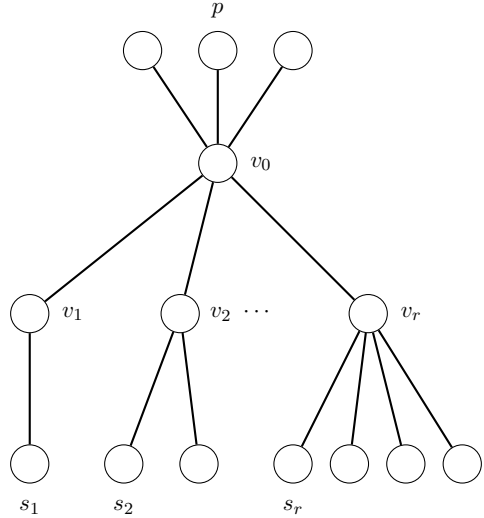
\end{centering}
\begin{Thr} \label{diam4T}
Conjecture \ref{TheConj} holds for trees of diameter 4.
\end{Thr}
\noindent
\textbf{Proof:}
If $T$ has diameter $4$, then $n \geq 5$.
Consider the prototypical
diameter $4$ tree
shown in Figure~\ref{examplefig7}.
We may assume there is a root $v_0$,
having $p \geq 0$ neighbors that are leaves, and $r \geq 2$ intermediate neighbors
$v_i$, $1 \leq i \leq r$,
each adjacent to $s_i \geq 1$ leaves.
The number of vertices of $T$ is
$n = p+1+r+\sum_{i=1}^r s_i$.
The number of leaves is $p+\sum_{i=1}^r s_i$.
Now consider applying
{\tt Algorithm LDiagonalize}$(T,x)$ with $x=-2+\frac{2}{n}$.
At initialization, each vertex is assigned
$d(v) = \deg(v) -2 + \frac{2}{n}$.
In particular, each leaf receives the permanent value $-1+\frac{2}{n} < 0$.
If $p > 0$ or $s_i > 1$,
for some $i$, then the number of
leaves (and hence the number of negative values)
is at least $\left\lceil \frac{n}{2} \right\rceil$.
Therefore we may assume $p=0$ and $s_i=1$ for all $i=1,\ldots,r$.
Applying the algorithm, we see that for all $i=1,\ldots,r$
the final value of $d(v_{i})$ is
$$
\deg(v_i) - 2 +\frac{2}{n} - \frac{1}{-1+\frac{2}{n}}
= \frac{2}{n} + \frac{1}{1-\frac{2}{n}} =\frac {{n}^{2}+2\,n-4}{n \left( n-2 \right) },
$$
which is strictly positive.
It suffices to show the final value of $d(v_0)$ is negative.
We have
\begin{eqnarray*}
d(v_0)& = &r - 2 + \frac{2}{n} + \sum_{i=1}^r \frac{1}{d(v_i)}
       = r - 2 + \frac{2}{n} - r\left(\frac {n \left( n-2 \right) }{{n}^{2}+2\,n-4}\right) \\
& =&  r - 2 + \frac{2}{n} - r\left(1 + \frac{ -4n+4}{{n}^{2}+2\,n-4}\right)
        = -2 +\frac{2}{n} + r\frac{4n-4}{{n}^{2}+2\,n-4}.
\end{eqnarray*}
Using $n=2r+1$, the last expression can be written as
$$
-2\,{\frac { \left( n-1 \right)  \left( 3\,n-4 \right) }{n \left( {n}^ {2}+2\,n-4 \right) }}
$$
which is clearly negative, completing the proof.
$\Box$

We see in the above proof that when the number of leaves is at least as large
as the number of internal vertices, then the negative values on the leaves are sufficient
to guarantee the inequality in the conjecture.
\begin{Cor}\label{Lotsofleaves}
If $T$ is  a tree of order $n$ having $p$ leaves,
and $p \geq \theratio{n}{2}$, then Conjecture \ref{TheConj} holds for $T$.
\end{Cor}

The website given in \cite{VCDT}
is a front-end to a database designed to give spectral information on trees, involving
classical, Laplacian, and normalized Laplacian eigenvalues.
Its implementation utilized B. McKay's list of trees \cite{McK},
and the high-performance Linear Algebra software LAPACK (C version) \cite{Lapack}.
The website is a nice tool for understanding eigenvalue distribution of trees
and formulating conjectures.
In fact, Conjecture~\ref{TheConj} was formulated using this tool,
and our experiments show that it holds for all trees
having order $n \leq 20$.

We observe that in general, trees with small diameter seem to have more Laplacian eigenvalues
below the average degree, and trees with large diameter seem to have fewer.
The extremal cases are the path $P_n$, having diameter $n-1$ and
exactly $\theratio{n}{2}$
below the average,
and the star $S_n$, which has diameter 2 and
has $n-1$ eigenvalues  in $[0,\overline{d})$.
For trees of diameter $d$, it is known \cite{Gro1990} that
$$
m_T(2,n]\geq \floorratio{d}{2}
$$
also providing evidence that trees with large diameter are more difficult.

Another indication that the problem might be hard
is that eigenvalues not smaller than $1$ abound.
In fact, by Theorem~\ref{thr-gamma} and $\gamma \leq \floorratio{n}{2}$,
one sees that $m_T[1,n] \geq \displayratio{n}{2}$.

On the positive side, it is shown in \cite{Braga2013}
that for any $T$ tree of order $n$
$$
m_T[0,2) \geq \displayratio{n}{2} .
$$
Since
$$
\lim_{n \rightarrow \infty} \overline{d}
= \lim_{n \rightarrow \infty}
2 -\frac{2}{n}
= 2
$$
for large trees, the conjecture is to within $\epsilon$ of being satisfied.
However, there do exist trees with eigenvalues between $\overline{d}$ and 2.

We finish by mentioning another open problem for trees.
Conjecture~\ref{TheConj}
was first stated in \cite{Tre2011} in the context of Laplacian energy.
We recall that the \emph{Laplacian energy} of a graph $G$ whose Laplacian spectrum is given
by $\mu_1 \geq \ldots  \geq \mu_n$, is defined \cite{Gut2006} as
$$
LE(G)=\sum_{i=1}^n |\mu_i -\overline{d}|,
$$
where $\overline{d} = \frac{2m}{n}$ is the average degree of $G$.
It is easy to show that
$$
LE(G)=2 \sum_{i=1}^{\sigma}\mu_i-2\sigma\overline{d},
$$
where $\sigma$ is the index of the smallest Laplacian eigenvalue greater than or equal to $\overline{d}$.
Hence, in order to compute or to bound the Laplacian energy of graphs,
one would like to know the number of eigenvalues above or below the average degree.
In \cite{Fri2011} it is shown that the tree of order $n$
having largest Laplacian energy is the star $S_n$.
It is not known which tree has the smallest Laplacian energy.
Some believe that the path $P_n$ has the smallest Laplacian energy among all trees of order $n$.
Settling Conjecture \ref{TheConj} could provide insight into
this problem.

\section{Final remarks}
In the absence of a solution
to the conjecture (affirming or providing a counterexample)
perhaps the strongest partial result would be to show the conjecture holds
for all but finitely many trees.
A weaker result would be that for each fixed $d$ among the trees of diameter $d$,
all but finitely many satisfy the conjecture.
Or, similarly, among trees with maximum degree $\Delta$,
the conjecture is satisfied for sufficiently large $n$.
Characterizing those trees for which $m_T[0,\bar{d}) = \theratio{n}{2}$
could also supply insight to the overall conjecture,
as well as determining how large $m_T (\bar{d},2)$ can be.
Finally, it is reasonable to seek a random spectral solution,
as other authors (e.g. \cite{Schwenk71,Bhamidi2012}) have done,
and show the conjecture holds for almost all trees.

\bibliographystyle{amsplain}
\bibliography{paper}

\end{document}